\documentclass[11pt]{article}
\usepackage{amssymb,amsmath}

\usepackage{latexsym}

\newfont{\bb}{msbm10 at 11pt}
\newfont{\bbsmall}{msbm8 at 8pt}
\usepackage[margin=1.4in]{geometry}

\newcommand{\wt}{\widetilde}
\newcommand{\R}{\mbox{\bb R}}

\newcommand{\C}{\mbox{\bb C}}

\newcommand{\N}{\mbox{\bb N}}

\newcommand{\B}{\mbox{\bb B}}

\newcommand{\esf}{\mbox{\bb S}}

\newcommand{\rth}{\R^3}

\newcommand{\tM}{\widetilde{M}}
\newcommand{\Mint}{M_{\infty}}

\newcommand{\qed}{\hfill \ensuremath{\Box}}

\def\S{{\Sigma}}
\def\s{{\sigma}}

\def\g{{\gamma}}

\def\l{{\lambda}}

\def\ve{{\varepsilon}}

\def\centerbmp#1#2#3{\vskip#2\relax\centerline{\hbox to#1{\special
    {bmp:#3 x=#1, y=#2}\hfil}}}

\newtheorem{theorem}{Theorem}[section]

\newtheorem{remark}[theorem]{Remark}

\newtheorem{definition}[theorem]{Definition}

\newenvironment{proof}{\smallskip\noindent{\it Proof.}\hskip \labelsep}
                          {\hfill\penalty10000\raisebox{-.09em}{$\Box$}\par\medskip}

\begin{document}
\begin{title}
{The rigidity of embedded constant mean curvature surfaces}
\end{title}
\begin{author}
{William H. Meeks, III\thanks{ This material is based upon
 work for the NSF under Award No. DMS -
 0405836. Any opinions, findings, and conclusions or recommendations
 expressed in this publication are those of the authors and do not
 necessarily reflect the views of the NSF.} \and Giuseppe Tinaglia}
\end{author}
\maketitle
\begin{abstract} We study the rigidity of complete, embedded
constant mean curvature surfaces in $\rth$. Among other things, we
prove that when such a surface has finite genus, then intrinsic
isometries of the surface extend to isometries of $\rth$ or its
isometry group contains an index two subgroup of isometries that
extend.

\vspace{.1cm} \noindent{\it Mathematics Subject Classification:}
Primary 53A10, Secondary 49Q05, 53C42

 \end{abstract}

\section{Introduction.}
In this paper we discuss some global results for certain complete
embedded surfaces $M$ in $\rth$ which have constant mean
curvature\footnote{We require that $M$ is equipped with a Riemannian
metric and that the inclusion map $i\colon M\to \rth$ preserves this
metric.}. If this mean curvature is zero, we call $M$ a minimal
surface and if it is nonzero, we call $M$ a $CMC$ surface. Our main
theorems deal with the rigidity of complete, embedded constant mean
curvature surfaces in $\rth$ with finite genus.

Recall that an isometric immersion $f\colon \S \to \rth$ of a
Riemannian surface  $\S$ is {\it congruent}  to another isometric
immersion $h\colon \S \to \rth$, if there exists an isometry
$I\colon \rth \to \rth$ such that $f = I\circ h$. We say that an
isometric immersion $f\colon \S\to \rth$ with constant mean
curvature $H$ is {\it rigid}, if whenever $h\colon \S\to \rth$ is
another isometric immersion with constant mean curvature $H$ or
$-H$, then $f$ is congruent to $h$.

In general, if $f\colon M \to \rth$ is an isometric immersion of a
simply-connected surface with constant mean curvature $H$ and $f(M)$
is not contained in a round sphere or a plane, then there exists a
smooth one-parameter deformation of the immersion $f$ through
non-congruent isometric immersions with mean curvature $H$; this
family contains all noncongruent isometric immersions of $M$ into
$\rth$ with constant mean curvature $H$ or $-H$. Thus, the rigidity
of simply-connected, constant mean curvature immersed surfaces fails
in a rather natural way. On the other hand, the main theorems
presented in this paper affirm that some complete,
nonsimply-connected, embedded constant mean curvature surfaces in
$\rth$ are rigid. More precisely, we prove the following theorems.

\begin{theorem}[Finite Genus Rigidity Theorem]\label{t8}
Suppose $M\subset \rth$ is a complete, embedded, constant mean
curvature surface of finite genus.
\begin{enumerate}
\item If $M$ is a minimal surface which is not the helicoid, then $M$ is rigid.
\item If $M$ is a $CMC$ surface with bounded Gaussian curvature, then $M$ is rigid.
\end{enumerate}
\end{theorem}

\begin{theorem}[Finite Genus Isometry Extension Theorem]\label{t1}
Let $M\subset \rth$ be a complete, embedded, constant mean curvature
surface of finite genus and let $\sigma\colon M\to M$ be an
isometry.
\begin{enumerate}
\item If $M$ has bounded Gaussian curvature or if $M$ is minimal, then $\sigma$ extends
to an isometry of $\rth$.
\item If $\sigma$ fails to extend to an isometry of $\rth$, then the
isometry group of $M$ contains a subgroup of index two,
consisting of those isometries which do extend to $\rth$. In particular, if
$\sigma$ fails to extend, then $\sigma^2$ extends.
\end{enumerate}
\end{theorem}

The first relevant result in the direction of revealing the rigidity
of certain constant mean curvature surfaces is a theorem of Choi,
Meeks, and White. In~\cite{cmw1} they proved that any properly
embedded minimal surface in $\rth$ with more than one end admits a
unique isometric minimal immersion into $\rth$. One of the
outstanding conjectures in this subject states that, except for the
helicoid, the inclusion map of a complete, embedded constant mean
curvature surface $M$ into $\rth$ is the unique such isometric
immersion with the same constant mean curvature up to congruence.
Since extrinsic isometries of the helicoid extend to ambient
isometries, the validity of this conjecture implies the closely
related conjecture that the intrinsic isometry group of any
complete, embedded constant mean curvature surface in $\rth$ is
equal to its ambient symmetry group. These two rigidity conjectures
were made by Meeks; see Conjecture 15.12 in \cite{mpe2} and the
related earlier Conjecture 22 in \cite{me14} for properly embedded
minimal surfaces. The theorems presented in this paper demonstrate
the validity of these rigidity conjectures under some additional
hypotheses.

The proofs of Theorems~\ref{t8} and \ref{t1} rely on the
classification of isometric immersions of simply-connected constant
mean curvature surfaces in $\rth$. A key ingredient in the proof of
these theorems is our Dynamics Theorem for $CMC$ surfaces in $\rth$
together with the Minimal Element Theorem which are proven in
\cite{mt4}. Among other things, in \cite{mt4} we prove that, under
certain hypotheses, a $CMC$ surface in $\rth$ contains an embedded
Delaunay surface\footnote{The Delaunay surfaces are $CMC$ surfaces
of revolution which were discovered and classified by
Delaunay~\cite{de1} in 1841. When these surfaces are embedded, they
are called {\it unduloids} and when they are nonembedded, they are
called {\it nodoids}.} at infinity. The fact that embedded Delaunay
surfaces are rigid is applied in the proofs of our main theorems.

Additionally, using techniques similar to those applied to prove
Theorems~\ref{t8} and \ref{t1}, we also demonstrate the following
related rigidity theorem. This theorem is well known in the special
case that the surface has finite topology (see for example,
\cite{ku2}).

\begin{theorem}\label{t2}
Suppose that $M\subset \rth$ is a complete embedded $CMC$ surface
with bounded Gaussian curvature.  If $$\liminf \left
(\frac{\mbox{\rm Area}[M\cap \B(R)]}{R^2} \right )=0 \quad \text{or}
\quad \liminf \left (\frac{\mbox{\rm Genus}[M\cap \B(R)]}{R^2}
\right )=0,$$ then $M$ is rigid.
\end{theorem}

\vspace{.1cm}

\noindent {\bf Acknowledgements} The authors would like to thank
Brian Smyth for helpful discussions.

\section{Background on the Dynamics Theorem.}\label{sec2}

Before stating results which we need from \cite{mt4} we introduce
some definitions; see also \cite{mpr10} for some related results for
minimal surfaces.

\begin{definition} \label{def} {\rm Suppose $W$ is a complete flat three-manifold
with boundary $\partial W=\S$ together with an isometric immersion
$f\colon W \to \rth$ such that $f$ restricted to the interior of $W$
is injective. We call the image surface $f(\S)$ a {\em strongly
Alexandrov embedded $CMC$ surface} if $f(\S)$ is a $CMC$ surface and
$W$ lies on the mean convex side of $\S$.}
\end{definition}

We note that, by elementary separation properties, any properly
embedded $CMC$ surface  in $\rth$ is always strongly Alexandrov
embedded. Furthermore, by item {\it 1} of Theorem~\ref{T} below, any
strongly Alexandrov embedded $CMC$ surface in $\rth$ with bounded
Gaussian curvature is properly immersed in $\rth$. We remind the
reader that the Gauss equation implies that a surface $M$ in $\rth$
with constant mean curvature has bounded Gaussian curvature if and
only if its principal curvatures are bounded in absolute value;
thus, $M$ having bounded Gaussian curvature is equivalent to $M$
having bounded second fundamental form.

\begin{definition} {\rm Suppose $M\subset \rth$ is a connected,
strongly Alexandrov embedded $CMC$
surface with bounded Gaussian curvature.
\begin{enumerate}
\item ${\cal T}(M)$ is the set of all connected, strongly
Alexandrov embedded $CMC$ surfaces $\S \subset \rth$, which are
obtained in the following way. There exists a sequence of points
$p_n\in M$, $\lim_{n\to \infty}|p_n|=\infty$, such that the
translated surfaces $M-p_n$ converge $C^2$ on compact subsets of
$\rth$ to a strongly Alexandrov embedded $CMC$ surface $\Sigma'$,
and $\Sigma$ is a connected component of $\Sigma'$ passing through
the origin. Actually we consider the immersed surfaces in ${\cal
T}(M)$ to be {\it pointed} in the sense that if such a surface is
not embedded at the origin, then we consider the surface to
represent two different surfaces in ${\cal T}(M)$ depending on a
choice of one of the two preimages of the origin.
\item $\Delta \subset  {\cal T}(M)$ is called {\em ${\cal T}$-invariant}, if $\S\in\Delta$
implies ${\cal T}(\S)\subset \Delta$.
\item A nonempty subset $\Delta\subset {\cal T}(M)$ is called a {\em minimal} ${\cal T}$-invariant
set, if it is ${\cal T}$-invariant and contains no smaller nonempty ${\cal
T}$-invariant subsets.
\item If $\S \in {\cal T}(M)$ and $\S$ lies in a minimal ${\cal T}$-invariant
subset of ${\cal T}(M)$, then $\S$ is called a  {\em minimal
element}
 of ${\cal T}(M)$.
\end{enumerate}}
\end{definition}

The following theorem is a collection of results and special cases
of statements taken from the $CMC$ Dynamics Theorem, the Minimal
Element Theorem and Theorem~4.1 in \cite{mt4}; in the statement of
this theorem, $\B(R)$ denotes the open ball of radius $R$ centered
at the origin in $\rth$.

\begin{theorem}\label{T}
Let $M$ be a connected, noncompact, strongly Alexandrov embedded
$CMC$ surface with bounded Gaussian curvature. Then:
\begin{enumerate}
\item $M$ is properly immersed in $\rth$.
\item ${\cal T}(M)$ is nonempty and ${\cal T}$-invariant.
\item Every nonempty ${\cal T}$-invariant subset of ${\cal T}(M)$ contains a nonempty minimal
${\cal T}$-invariant subset. In particular, since ${\cal T}(M)$ is
itself a nonempty ${\cal T}$-invariant set, ${\cal T}(M)$ always
contains minimal elements.
\item Suppose $$\liminf \left(\frac{\text{\rm Area}[M\cap \B(R)]}{R^2}\right)=0 \quad \text{or}
\quad \liminf \left(\frac{\text{\rm Genus}[M\cap \B(R)]}{R^2}\right)=0.$$
Then ${\cal T}(M)$ contains a minimal element $\Sigma$  which is an embedded Delaunay surface.
\item Suppose $M$ has a plane of Alexandrov symmetry and more than one end.
Then ${\cal T}(M)$ contains a minimal element $\Sigma$ which is an
embedded Delaunay surface.
\end{enumerate}
\end{theorem}

\section{Background on Calabi's and Lawson's Rigidity Theorems.}\label{sec2b}
In this section we review the classical rigidity theorems of Calabi
and Lawson for simply-connected constant mean curvature surfaces in
$\rth$ (see also~\cite{Bonnet}). The rigidity theorem of Lawson is motivated by the earlier
result of Calabi~\cite{ca1} who classified the set of isometric
minimal immersions of a simply-connected Riemannian surface $\Sigma$
into $\rth$; we now describe Calabi's classification theorem.

Suppose $f\colon \S \to \rth$ is a isometric minimal immersion and
$\S$ is simply-connected. In this case the coordinate functions
$f_1, f_2, f_3$ are harmonic functions which are the real parts of
corresponding holomorphic functions $h_1, h_2, h_3$ defined on $\S$.
For any $\theta \in [0,2\pi)$, the map $f_\theta= {\rm
Re}(e^{i\theta}(h_1,h_2,h_3))\colon \S \to \rth$ is an isometric
minimal immersion of $\S$ into $\rth$; the immersions $f_\theta$ are
called {\em associate immersions} to $f$. Many classical examples of
minimal surfaces arise from this associate family construction. For
example, simply-connected regions on a catenoid are the images of
regions in the helicoid under the associate map for
$\theta=\frac{\pi}{2}$; in this case the corresponding coordinate
functions on these domains are conjugate harmonic functions and
consequently, the catenoid and the helicoid are called {\em
conjugate} minimal surfaces. 

Calabi's classification or rigidity
theorem states that if $\Sigma$ is not flat, then for any isometric
minimal immersion $F\colon \S \to \rth$, there exists a {\em unique}
$\theta \in [0,\pi)$ such that $F$ is congruent to $f_\theta$, i.e.
there exists an isometry $I\colon \rth\to\rth$ such that as
mappings, $F=I\circ f_\theta$. This notion of rigidity does {\it
not} mean that the image surface $f(\S)$ cannot be congruent to the
image of an associate surface $f_\theta(\S)$, where $\theta \in
(0,\pi)$. For example, let $f\colon \C\to\rth$ be a parametrization
of the classical Enneper surface and let $f_{\frac{\pi}{2}}$ be the
conjugate mapping. Then the images of these immersions are congruent
as subsets of $\rth$ but these immersions are not congruent as
mappings. In fact, the rotation $R_{\frac{\pi}{2}}$ counter
clockwise by $\frac{\pi}{2}$ in the usual parameter coordinates $\C$
for Enneper's surface is an intrinsic isometry of this surface which
does not extend to an isometry of  $\rth$ and $f\circ
R_{\frac{\pi}{2}}$ is congruent to the immersion
$f_{\frac{\pi}{2}}$.

Calabi condition for rigidity of $\S$ is equivalent to the property
that the conjugate harmonic coordinate functions are well defined,
which by Cauchy's theorem is equivalent to the property that the
flux vector\footnote{The flux of an oriented unit speed curve
$\gamma\subset \S$ is the integral of $J(\gamma')$ along $\gamma$,
where $J$ is the complex structure.} of any simple closed curve on
$\Sigma$ is zero.

Lawson's Rigidity Theorem that we referred to in the previous
paragraph appears in Theorem 8 in \cite{la3} and holds in space
forms other than $\rth$ as well. We will not need his theorem in its
full generality and we state below the special case which we will
apply in the next section.

\begin{theorem}[Lawson's Rigidity Theorem  for $CMC$ surfaces in
$\rth$]\label{Law} If $f\colon \S \to \rth$ is an isometric $CMC$
immersion with mean curvature $H$ and $\Sigma$ is simply-connected,
then there exists a differentiable $2\pi$-periodic family of
isometric immersions
$$f_\theta \colon \S\to\rth$$ of constant mean curvature $H$ called
associate immersion to $f$. Moreover, up to congruences the maps
$f_\theta$, for $\theta\in[0,2\pi)$, represent all isometric
immersions of $\S$ into $\rth$ with constant mean curvature $H$ or $-H$ and
these immersions are non-congruent to each other if $f(M)$ is not
contained in a sphere.
\end{theorem}

Note that if $A$ represents the second fundamental form of $f$ and
$A_\theta$ the one of $f_\theta$, these forms are related by the
following equation, see \cite{smy}:
\begin{equation}\label{atheta}
A_\theta= \cos (\theta)(A-HI)+\sin(\theta)J(A-HI)+HI,
\end{equation} where $I$ is the identity matrix and $J$ is the
almost complex structure on $M$.

\section{Proofs of the main theorems.}\label{sec3}
Our rigidity theorems are motivated by several classical results on
the rigidity of certain complete embedded minimal surfaces. The
first relevant result in this direction is a theorem of Choi, Meeks,
and White~\cite{cmw1} who proved that any properly embedded minimal
surface in $\rth$ with more than one end admits a unique isometric
minimal immersion into $\rth$; their result proved a special case of
the conjecture of Meeks  \cite{me14} that the inclusion map of a
properly embedded, nonsimply-connected minimal surfaces in $\rth$ is
the unique minimal immersion of the surface into $\rth$ up to
congruence.

We now prove the following important special case of Theorem~\ref{t8}
regarding properly embedded minimal surfaces.%

\begin{theorem}\label{t3} If $M$ is a connected,
properly embedded minimal surface in $\rth$ with finite genus which is not a
helicoid, then $M$ is rigid.
\end{theorem}
\begin{proof}
If $M$ has more than one end, then the result of Choi, Meeks, and
White implies that $M$ satisfies the conclusions of the theorem.

Suppose $M$ has finite genus, one end and $M$ is not a helicoid. The
main theorem of Meeks and Rosenberg in \cite{mr8} then implies that
$M$ is a plane or is asymptotic to a helicoid $\S$. In the latter
case, $M$ can be thought of as being a helicoid with a finite
positive number of handles attached\footnote{This asymptotic to a
helicoid property of $M$ follows directly from the proof of the
uniqueness of the plane and the helicoid as being the only properly
embedded simply-connected minimal surfaces. The proof of this
asymptotic argument is explained briefly in the last section of
\cite{mr8}. In the last section of their paper, Meeks and Rosenberg
inadvertently describe a stronger analytic statement which is false.
In any case, a simplification of the asymptotic to a helicoid
property of $M$ is also proved in \cite{mpe3}}. Since
Theorem~\ref{t3} holds for planes, assume now that $M$ is not a
plane. Since $M$ is asymptotic to a helicoid $\S$, any plane $P$
orthogonal to the axis of $\S$ intersects $M$ in an analytic set
with each component of $M\cap P$ having dimension one and such that
outside of some ball in $\rth$, $M\cap P$ consists of two proper
arcs asymptotic to the line $\S\cap P$. Since $M$ has finite
positive genus, elementary Morse theory implies that for a certain
choice of $P$, $M\cap P$ is a one-dimensional analytic set with a
vertex contained in the intersection of two open analytic arcs in
$M\cap P$. An elementary combinatorial argument implies that $P\cap
M$ contains a simple closed oriented curve $\Gamma$ bounding a
compact disk whose interior is disjoint from $M$. It follows that
the integral of the conormal to $\Gamma$ has a nonzero dot product
with the normal to $P$. The existence of $\Gamma$ implies that for
at least one of the coordinate functions $x_i$ of $M$, the conjugate
harmonic function of $x_i$ is not well-defined (for example, if
$\Gamma$ lies in the $(x_1,x_2)$-plane, then the conjugate harmonic
function of the $x_3$-coordinate function is not well-defined on
$M$). From our discussion of the Calabi Rigidity Theorem in the
previous section, it follows that the inclusion map of $M$ into
$\rth$ is the unique isometric minimal immersion of $M$ into $\rth$
up to congruence. This completes the proof of the theorem.
\end{proof}

We will now apply the above theorem and the results described in
Sections~\ref{sec2} and \ref{sec2b} to prove Theorems~\ref{t8},
\ref{t1} and \ref{t2}.

\vspace{.2cm}

In what follows, let $M\subset \rth$ be a complete, embedded,
constant mean curvature surface of finite genus.

\vspace{.3cm}

\noindent{\bf Proof of Theorem~\ref{t8}.} We first prove item {\it
1}. Suppose $M$ is a minimal surface which is not a helicoid. If $M$
has positive injectivity radius, then $M$ is properly embedded by
the Minimal Lamination Closure Theorem in \cite{mr13} and the result
is a consequence of Theorem \ref{t3}. If $M$ fails to have positive
injectivity radius, then the local picture theorem on the scale of
the topology in \cite{mpr14} implies that there exists a sequence of
compact domains $\Delta_n\subset M$ such that, after scaling and a
rigid motion of $\rth$, the new domains $\wt{\Delta}_n$ converge
smoothly with multiplicity one to a properly embedded genus zero
minimal surface $M_\infty$ in $\rth$ with injectivity radius one
(and hence not simply-connected) or there exist closed geodesics
$\gamma_n \subset \wt{\Delta}_n$ with nontrivial flux (the integral
of the conormal is nonzero). The latter happens when the
$\wt{\Delta}_n$ converge smoothly away from two vertical lines $L_1$
and $L_2$ to a foliation ${\cal F}$ by horizontal planes of $\rth$.
These lines are the singular sets of convergence to ${\cal F}$. Such
a picture is called a parking garage structure on $\rth$ and the
geodesics $\gamma_n$ correspond to ``connection'' curves between the
``columns'' $L_1$ and $L_2$ (see \cite{mpr14}). On the one hand, if
$M_\infty$ is a properly embedded genus zero minimal surface with
injectivity radius one, then by Theorem 4.1, $M_\infty$ is rigid and
a compactness argument implies that $M$ is also rigid (see the proof
of Theorem \ref{t2} for this type of argument). On the other hand,
minimal surfaces with nontrivial flux are rigid, by our discussion
in Section 3. This completes the proof of item {\it 1}.

Since finite genus implies that $\liminf (\frac{\text{Genus}[M\cap
\mathbb{B}(R)]}{R^2})=0$, item {\it 2}\, is a consequence of
Theorem~\ref{t2}, which is proved below. This completes the proof of
Theorem~\ref{t8}. \qed

\vspace{.3cm}

\noindent{\bf Proof of Theorem~\ref{t1}.} Since rigidity implies
that intrinsic isometries extend to extrinsic isometries, item {\it
1} is a simple consequence of Theorem \ref{t8}. We now prove item
{\it 2}. Suppose $M$ is a $CMC$ surface and {\bf does not} have
bounded Gaussian curvature. Let $i\colon M \to \rth$ be the
inclusion map.  If $\sigma$ is an isometry of $M$, then $i\circ
\sigma$ is congruent to a unique associate surface
$i_{\theta(\sigma)}$,  $\theta(\sigma)\in [0,2\pi)$. By defining
$F(\sigma)=\theta(\sigma)$, we obtain a homomorphism $F \colon
\text{\rm Isom}(M) \to \esf^1$. We are going to show that $F \colon
\text{\rm Isom}(M)\subset \{0, \pi \} = \mathbb{Z}_2\subset \esf^1$
which will imply item {\it 2} of the theorem.

Assume that $\s\colon M\to M$ is an isometry that does not extend to
$\rth$. By Theorem \ref{Law} and the previous discussion we know
that $i\circ\s\colon M\to \rth$ must be congruent to an associate
immersion $i_{\theta}\colon M\to \rth$, where $\theta\in (0, 2\pi)$.

The local picture
theorem on the scale of curvature in
\cite{mpr10} states that there exists a sequence of points $p_n\in
M$ and positive numbers $\ve_n,\l_n$ such that:

\begin{enumerate}
\item $\lim_{n\to\infty}{\ve_n}\to 0$,
$\lim_{n\to\infty}{\l_n}=\infty$ and
$\lim_{n\to\infty}{\l_n\ve_n}=\infty$.
\item The component $M_n$ of $M\cap B(p_n,\ve_n))$ that
contains $p_n$ is compact with $\partial M_n\subset\partial
B(p_n,\ve_n)$. \item The second fundamental forms of the surfaces
$\tM_n=\l_n M_n\subset \l_n B(p_n,\ve_n)\subset \rth$ are uniformly
bounded and are equal to one at the related points $\wt{p}_n$.
\item The translated surfaces $\tM-\wt{p}_n$ converge
with multiplicity one to a connected, properly embedded minimal
surface $\Mint \subset\rth$ with bounded curvature and genus zero.
\end{enumerate}

Suppose first that $\Mint$ is a helicoid. We will use the
embeddedness property of $M$ to show that $\theta=\pi$.  If $\Mint$
is a helicoid, then the associate surfaces $(\tM-\wt{p}_n)_{\theta}$
can be chosen to approximate a large region of the related associate
surface $(\Mint)_{\theta}$ to the helicoid $\Mint$; this can be seen
by letting $H$ go to zero in equation \eqref{atheta}. If $\theta\neq
\pi$, then $(\Mint)_{\theta}$ intersects itself transversely which
means that $(\tM-\wt{p}_n)_{\theta}$ is not embedded. This
contradiction concludes the argument.

It remains to show the case when $\Mint$ is not the helicoid.  If
$\Mint$ is not a helicoid, then Theorem~\ref{t3} implies
$({\Mint})_{\theta}$ is not well-defined unless $\theta = \pi$.
\qed

\vspace{.3cm}

\noindent{\bf Proof of Theorem~\ref{t2}.} Suppose $M$ is a $CMC$
surface with bounded Gaussian curvature such that either $\liminf
(\frac{\text{Area}[M\cap \mathbb{B}(R)]}{R^2})=0$ or $\liminf
(\frac{\text{Genus}[M\cap \mathbb{B}(R)]}{R^2})=0$. Without loss of
generality we will assume $H_M=1$. By item {\it 3} of
Theorem~\ref{T}, ${\cal T}(M)$ contains an embedded Delaunay surface
$\Sigma$. More precisely, for $n\in \N$, there exist compact annular
domains $\Delta_n\subset M$ and points $p_n\in \Delta_n$ such that
the translated surfaces $\Delta_n-p_n$ converge $C^2$ to $\Sigma$ on
compact subsets of $\rth$. For concreteness, suppose $g_1$ denotes
the inclusion map of $\S$ into $\rth$. First we show that the
immersion $g_1$ is rigid.

Let $\pi\colon\wt{\S}\to\S$ denote the universal covering of $\S$.
Consider $\wt{\S}$ with the induced Riemannian metric and let
$f=g_1\circ\pi\colon \wt{\S}\to \rth$ be the related isometric
immersion. Let $f_{\theta}\colon \wt{\S}\to \rth$ be the associate
immersion for angle $\theta \in [0,2\pi)$, given in
Theorem~\ref{Law}; note $f_0=f$. Suppose $g_2$ is another isometric
immersion of $\S$ into $\rth$ which is not congruent to $g_1$. This
implies that $g_2\circ \pi$ is congruent to $f_{\overline{\theta}}$
for a certain $\overline{\theta} \in (0,2\pi)$.

Let $\widetilde{\gamma}\subset \wt{\S}$ be a lift\footnote{The curve
$\wt{\g}$ is a compact embedded arc in $\wt{\S}$ which is the image
of a lift of map $\g\colon [0,1]\to \S$.} of the shortest geodesic circle
$\g\subset\Sigma$. We will prove that for any $\theta\in(0,2\pi)$ the
endpoints of $f_\theta(\widetilde\gamma)$ are distinct, which means
that the associate immersions to $g_1$ do not exist. We will
accomplish this by describing the geometry of
$f_\theta(\widetilde{\gamma}).$

A computation shows that for the geodesic $\widetilde{\gamma}$ and
the immersion $f_\theta$, the curvature $k_\theta$ and torsion $\tau_\theta$
of $f_\theta(\widetilde{\gamma})(t)$ are given by

\begin{equation}\label{formula}
k_\theta=\langle A_\theta (\widetilde{\gamma}'(t)),\widetilde{\gamma}'(t)\rangle
\text{ and } \tau_\theta=-\langle A_\theta
(\widetilde{\gamma}'(t)),J(\widetilde{\gamma}'(t))\rangle.
\end{equation}
Furthermore, since $\widetilde{\gamma}$ is a lift of the shortest
geodesic circle in $\S$, there exist $s\leq 0<2\leq r$, $r+s=2$ such
that in the $\widetilde{\gamma}'$, $J\widetilde{\gamma}'$ basis, the
second fundamental form along $\widetilde{\gamma}$ is expressed as
the matrix
$$A=\left( \begin{array}{cc}
r & 0\\
0 & s
\end{array}\right).$$

Consequently, equations \eqref{atheta} and \eqref{formula} imply
that
\begin{equation}
k_\theta=\cos(\theta)(r-1)+1 \;\text{ and } \; \tau_\theta=\sin(\theta)(1-r).
\end{equation}
In particular, $k_\theta$ and $\tau_\theta$ are constants. If
$\theta\neq \pi$, then $f_\theta(\widetilde{\gamma})$ is contained
in a helix, while if $\theta=\pi$, it is contained in circle of
radius $R=|2-r|=|s|<|r|$. Since the length of
$f_\pi(\widetilde{\gamma})$ is $\frac{2\pi}{r}$, it follows that in
either case the endpoints of $f_\theta(\widetilde{\gamma})$ are
distinct.

Since the compact annuli $\Delta_n-p_n$ converge $C^2$ to the
embedded Delaunay surface $\S$ as $n\to\infty$, we conclude that the
associate immersions for $\theta\in(0,2\pi)$ are not well-defined on
$\Delta_n-p_n$ for $n$ large. By Theorem~\ref{Law}, Theorem~\ref{t8}
now follows.  \qed

\begin{remark} In the above proof, we showed that any embedded
Delaunay surface is rigid. However, the same computations prove that
if $f$ represents the inclusion map of a nodoid into $\rth$, then
the associate immersions $f_\theta$, $\theta \in (0,\pi)$ are never
well-defined and the associate immersions $f_\pi$ are well-defined
for an infinite countable set of nodoids.
\end{remark}

\begin{remark} The conclusions of the Finite Genus Rigidity Theorem
(Theorem~\ref{t8}) should hold without the hypothesis that $M$ have
bounded Gaussian curvature. This improvement would be based on
techniques we are developing in \cite{mt1} to prove curvature
estimates for certain embedded $CMC$ surfaces in $\rth$.
\end{remark}

\begin{remark}
The full generality of the $CMC$ Dynamics Theorem and the Minimal
Element Theorem can be used to prove rigidity of embedded $CMC$
surfaces under hypotheses which imply the existence of an embedded
Delaunay surface at ``infinity.'' For instance, item {\it 4} of
Theorem \ref{T} implies that  if $M$ has bounded Gaussian curvature,
a plane of Alexandrov symmetry and more than one end, then $M$ is
rigid.
\end{remark}

\begin{remark}
It also makes sense to talk about ``local rigidity'' of a $CMC$
surface . The full generality of the $CMC$ Dynamics Theorem and the
Minimal Element Theorem can be used to prove local rigidity of
certain embedded $CMC$ surfaces (see~\cite{smyt1}).\end{remark}

\center{William H. Meeks, III at bill@math.umass.edu\\
Mathematics Department, University of Massachusetts, Amherst, MA,
01003}
\center{Giuseppe Tinaglia gtinagli@nd.edu   \\
Mathematics Department, University of Notre Dame, Notre Dame, IN,
46556-4618}

\bibliographystyle{plain}
\bibliography{tempbill}

\end{document}